\input amstex\documentstyle{amsppt}  
\pagewidth{12.5cm}\pageheight{19cm}\magnification\magstep1
\topmatter
\title From conjugacy classes in the Weyl group to unipotent classes, II\endtitle
\rightheadtext{Conjugacy classes in the Weyl group and unipotent classes, II}
\author G. Lusztig\endauthor
\address{Department of Mathematics, M.I.T., Cambridge, MA 02139}\endaddress
\thanks{Supported in part by the National Science Foundation}\endthanks
\endtopmatter   
\document
\define\hWW{\hat\WW}

\define\uuG{\un{\un G}}

\define\uWW{\un\WW}

\define\pe{\perp}
\define\si{\sim}

\define\sqc{\sqcup}

\define\qua{\quad}

\define\op{\oplus}
   
\redefine\sp{\spadesuit}
\define\part{\partial}

\define\m{\mapsto}
\define\do{\dots}

\define\sub{\subset}    

\define\T{\times}
\define\ti{\tilde}
\define\nl{\newline}
\redefine\i{^{-1}}

\define\un{\underline}

\define\a{\alpha}
\redefine\b{\beta}
\redefine\c{\chi}
\define\g{\gamma}

\define\e{\epsilon}
\define\et{\eta}
\define\io{\iota}

\define\p{\pi}

\define\ps{\psi}
\define\r{\rho}
\define\s{\sigma}
\redefine\t{\tau}
\define\th{\theta}
\define\k{\kappa}
\redefine\l{\lambda}
\define\z{\zeta}
\define\x{\xi}

\redefine\G{\Gamma}

\define\Ph{\Phi}
\define\Ps{\Psi}

\define\kk{\bold k}

\define\nn{\bold n}

\define\CC{\bold C}

\define\NN{\bold N}

\define\QQ{\bold Q}

\define\WW{\bold W}
\define\ZZ{\bold Z}

\define\ca{\Cal A}

\define\cp{\Cal P}
\define\cq{\Cal Q}
\define\car{\Cal R}
\define\cs{\Cal S}
\define\ct{\Cal T}

\define\cv{\Cal V}

\define\te{\ti e}

\define\tA{\ti A}

\define\tcp{\ti{\cp}}

\define\CA{Ca}
\define\KL{KL}
\define\LU{L1}
\define\ICC{L2}
\define\UNI{L3}
\define\USCI{L4}
\define\USCII{L5}
\define\USCIII{L6}
\define\HK{L7}
\define\WE{L8}
\define\CSM{L9}
\define\SPA{S1}
\define\SPAE{S2}
\define\SPAC{S3}
\define\XUE{Xue}

\define\uuL{\un{\un L}}
\define\X{\Xi}
\define\bV{\bar V}

\head Introduction\endhead
\subhead 0.1\endsubhead
Let $G$ be a connected reductive algebraic group over an algebraically closed field $\kk$ of characteristic 
$p\ge0$. Let $\uuG$ be the set of unipotent conjugacy classes in $G$. Let $\uWW$ be the set of conjugacy classes 
in the Weyl group $\WW$ of $G$. Let $\Ph:\uWW@>>>\uuG$ be the (surjective) map defined in \cite{\WE}. For 
$C\in\uWW$ we denote by $m_C$ the dimension of the fixed point space of $w:\cv@>>>\cv$ where $w\in C$ and $\cv$ 
is the reflection representation of the Coxeter group $\WW$. The following result provides a one sided inverse 
for $\Ph$.

\proclaim{Theorem 0.2} For any $\g\in\uuG$ the function $\Ph\i(\g)@>>>\NN$, $C\m m_C$ reaches its minimum at a 
unique element $C_0\in\Ph\i(\g)$. Thus we have a well defined map $\Ps:\uuG@>>>\uWW$, $\g\m C_0$ such that 
$\Ph\Ps:\uuG@>>>\uuG$ is the identity map.
\endproclaim
It is likely that when $\kk=\CC$, the map $\Ps$ coincides with the map $\uuG@>>>\uWW$ defined in 
\cite{\KL, Section 9}, which we will denote here by $\Ps'_0$ (note that $\Ps'_0$ has not been computed explicitly
in all cases). It is enough to prove the theorem in the case where 
$G$ is almost simple; moreover in that case it is enough to consider one group in each isogeny class. When $G$ 
has type $A$, the theorem is trivial ($\Ph$ is a bijection). For the remaining types the proof is given in 
Sections 1,2.

\subhead 0.3\endsubhead
Let $G_0$ be a connected reductive algebraic group over $\CC$ of the same type as $G$; we identify the Weyl group
of $G_0$ with $\WW$. Let $\uuG_0$ be the set of unipotent classes of $G_0$. Let $\Ph_0:\uWW@>>>\uuG_0$,
$\Ps_0:\uuG_0@>>>\uWW$ be the maps defined like $\Ph,\Ps$ (for $G_0$ instead of $G$). 

\proclaim{Theorem 0.4}(a) There is a unique (necessarily surjective) map $\r:\uuG@>>>\uuG_0$ such that 
$\Ph_0=\r\Ph$. We have $\r=\Ph_0\Ps$.

(b) There is a unique (necessarily injective) map $\p:\uuG_0@>>>\uuG$ such that $\Ps_0=\Ps\p$. We have 
$\p=\Ph\Ps_0$.
\endproclaim
The second sentence in 0.4(a) (resp. 0.4(b)) follows from the first since $\Ph\Ps=1$. The uniqueness of $\r$ 
(resp. $\p$) follows from the surjectivity of $\Ph$ (resp. injectivity of $\Ps$).  The surjectivity of $\r$ (if 
it exists) follows from the surjectivity of $\Ph_0$. The injectivity of $\p$ (if it exists) follows from the 
injectivity of $\Ps_0$. The existence of $\r$ (resp. $\p$), which is equivalent to the identity 
$\Ph_0=\Ph_0\Ps\Ph$ (resp. $\Ps_0=\Ps\Ph\Ps_0$) is proved in Section 3, where various other characterizations of 
$\r$ and $\p$ are given.

\subhead 0.5\endsubhead
Let $\hWW$ be the set of irreducible representations of $\WW$ (over $\QQ$) up to isomorphism. Let $\cs_\WW$ be 
the subset of $\hWW$ introduced in \cite{\LU}; it consists of representations of $\WW$ which were later called
"special representations". For any $E\in\hWW$ let $[E]$ be the unique object of $\cs_\WW$ such that $E,[E]$ 
are in the same two-sided cell of $\WW$; thus $E\m[E]$ is a surjective map $\hWW@>>>\cs_\WW$. For any 
$\g\in\uuG_0$ let $E_\g\in\hWW$ the Springer representation corresponding to $\g$. We define
$\ti\Ph_0:\uWW@>>>\cs_\WW$ by $C\m[E_{\Ph_0(C)}]$. Let $\uuG_0^{\sp}=\{\g\in\uuG_0;E_\g\in\cs_\WW\}$. It is 
known that $\g\m E_\g$ is a bijection $\uuG_0^{\sp}@>\si>>\cs_\WW$. This, together with the surjectivity of 
$\Ph_0$, shows that $\ti\Ph_0$ is surjective. We now define $\ti\Ps_0:\cs_\WW@>>>\uWW$ by 
$\ti\Ps_0(E_\g)=\Ps_0(\g)$, $\g\in\uuG_0^{\sp}$. Note that $\ti\Ps_0$ is injective and $\ti\Ph_0\ti\Ps_0=1$. We 
have the following result.

\proclaim{Theorem 0.6} (a) The (surjective) map $\ti\Ph_0:\uWW@>>>\cs_\WW$ depends only on the Coxeter group 
structure of $\WW$. In particular, $\ti\Ph_0$ for $G_0$ of type $B_n$ coincides with $\ti\Ph_0$ for $G_0$ of type
$C_n$ ($n\ge3$).

(b) The (injective) map $\ti\Ps_0:\cs_\WW@>>>\uWW$ depends only on the Coxeter group structure of $\WW$. In 
particular, $\ti\Ps_0$ for $G_0$ of type $B_n$ coincides with $\ti\Ps_0$ for $G_0$ of type $C_n$ ($n\ge3$).
\endproclaim
The proof is given in Section 4. Note that (a) was conjectured in \cite{\HK, 1.4}. 

Let $\uWW_{\sp}$ be the image of $\ti\Ps_0$ (a subset of $\uWW$). We say that $\uWW_{\sp}$ is the set of {\it 
special conjugacy classes} of $\WW$. Note that $\ti\Ph_0$ defines a bijection $\uWW_{\sp}@>\si>>\cs_\WW$. Thus, 
the set of special conjugacy classes in $\WW$ is in natural bijection with the set of special representations of 
$\WW$. From Theorem 0.6 we see that there is a natural retraction $\uWW@>>>\uWW_{\sp}$, $C\m\ti\Ps_0(\ti\Ph_0(C))$
which depends only on the Coxeter group structure of $\WW$.

\subhead 0.7\endsubhead
The paper is organized as follows. In Section 1 (resp. Section 2) we describe explicitly the map $\Ph$ in the 
case where $G$ is almost simple of classical (resp. exceptional) type and prove Theorem 0.2. In Section 3 we 
prove Theorem 0.4. In Section 4 we prove Theorem 0.6. In Section 5 we describe explicitly for each simple type 
the bijection $\uWW_{\sp}@>\si>>\cs_\WW$ defined by $\ti\Ph_0$.

\subhead 0.8. Notation\endsubhead
For $\g,\g'\in\uuG$ (or $\g,\g'\in\uuG_0$) we write $\g\le\g'$ if $\g$ is contained in the closure of $\g'$.

\head 1. Isometry groups\endhead
\subhead 1.1\endsubhead
Let $\cp^1$ be the set of sequences $p_*=(p_1\ge p_2\ge\do\ge p_\s)$ in $\ZZ_{>0}$. For $p_*\in\cp^1$ we set 
$|p_*|=p_1+p_2+\do+p_\s$, $\t_{p_*}=\s$, $\mu_j(p_*)=|\{k\in[1,\s];p_k=j\}|$ ($j\in\ZZ_{>0}$). Let 
$\cp^0=\{p_*\in\cp^1;\t_{p_*}=\text{ even}\}$. For $N\in\NN,\k\in\{0,1\}$ let 
$\cp^\k_N=\{p_*\in\cp^\k;|p_*|=N\}$. Let $\tcp=\{p_*\in\cp^0;p_1=p_2,p_3=p_4,\do\}$. Let $\cs^1$ be the set of 
$r_*\in\cp^1$ such that $r_k$ is even for all $k$. Let $\cs^0=\cs^1\cap\cp^0$. For $N\in\NN,\k\in\{0,1\}$ let 
$\cs^\k_N=\cs^\k\cap\cp^\k_N$.

We fix $\nn\in\NN$ and we define $n\in\NN$, $\k\in\{0,1\}$ by $\nn=2n+\k$. We set 
$\ca^k_{2n}=\{(r_*,p_*)\in\cs^\k\T\tcp;|r_*|+|p_*|=2n\}$.

\subhead 1.2\endsubhead
In the remainder of this section we fix $V$, a $\kk$-vector space of dimension $\nn=2n+\k$ (see 1.1) with a fixed
bilinear form $(,):V\T V@>>>\kk$ and a fixed quadratic form $Q:V@>>>\kk$ such that (i) or (ii) below holds:

(i) $Q=0$, $(x,x)=0$ for all $x\in V$, $V^\pe=0$;

(ii) $Q\ne0$, $(x,y)=Q(x+y)-Q(x)-Q(y)$ for $x,y\in V$, $Q:V^\pe@>>>\kk$ is injective.
\nl
Here, for any subspace $V'$ of $V$ we set $V'{}^\pe=\{x\in V;(x,V')=0\}$. In case (ii) it follows that $V^\pe=0$ 
unless $\k=1$ and $p=2$ in which case $\dim V^\pe=1$. Let $Is(V)$ be the subgroup of $GL(V)$ consisting of all
$g\in GL(V)$ such that $(gx,gy)=(x,y)$ for all $x,y\in V$ and $Q(gx)=Q(x)$ for all $x\in V$. In this section we 
assume that $G$ is the identity component of $Is(V)$.

\subhead 1.3\endsubhead
In this subsection we assume that $\nn\ge3$. Let $W$ be the group of permutations of $[1,\nn]$ that commute with 
the involution $\c:i\m\nn-i+1$. If $Q=0$ or if $\k=1$ we identify (as in \cite{\WE, 1.4, 1.5}) $\WW$ with $W$. If
$Q\ne0$ and $\k=0$ we identify (as in \cite{\WE, 1.4, 1.5}) $\WW$ with the group $W'$ of even permutations of 
$[1,\nn]$ commuting with $\c$; in this case let $\uWW_0$ (resp. $\uWW_1$) be the set of conjugacy classes in 
$\WW$ which are not conjugacy classes of $W$ (resp. form a single conjugacy class in $W$) and we denote by 
$\ti{\uWW}$ the set of $W$-conjugay classes in $W'$.

If $Q=0$ we identify $\uWW$ with $\ca^1_{2n}$ by associating to the 
conjugacy class of $w\in W$ the pair $(r_*,p_*)$ where $r_*$ is the multiset consisting of the sizes of cycles of
$w$ which are $\c$-stable and $p_*$ is the multiset consisting of the sizes of the cycles of $w$ which are not 
$\c$-stable. If $Q\ne0$ we identify $\uWW=\ca^1_{2n}$ (if $\k=1$) and $\ti{\uWW}=\ca^0_{2n}$ (if $\k=0$) by 
associating to the conjugacy class of $w$ in $W$ the pair $(r_*,p_*)$ where $r_*$ is the multiset consisting of 
the sizes of cycles of $w$ (other than fixed points) which are $\c$-stable and $p_*$ is the multiset consisting 
of the sizes of cycles of $w$ which are not $\c$-stable.

\subhead 1.4\endsubhead
Let $\ct_{2n}$ be the set of all $c_*\in\cp^1_{2n}$ such that $\mu_j(c_*)$ is even for any odd $j$.

Let $\ct^{(2)}_{2n}$ be the set of all pairs $(c_*,\e)$ where $c_*\in\ct_{2n}$ and $j\m\e(j)\in\{0,1\}$ is a 
function defined on the set $\{j\in\{2,4,6,\do\};\mu_j(c_*)\in\{2,4,6,\do\}\}$.

Let $\ti{\ct}^{(2)}_{2n}=\{(c_*,\e)\in\ct^{(2)}_{2n};\t_{c_*}\text{ is even}\}$.

Let $\cq$ be the set of all $c_*\in\cp^1$ such that  $\mu_j(c_*)$ is even for any even $j$. For $N\in\NN$ let 
$\cq_N=\cq\cap\cp^1_N$.

If $Q=0,p\ne2$ we identify $\uuG=\ct_{2n}$ by associating to $\g\in\uuG$ the multiset consisting of the sizes of 
the Jordan blocks of an element of $\g$. 

If $Q=0,p=2$ we identify $\uuG=\ct_{2n}^{(2)}$ by associating to $\g\in\uuG$ the pair $(c_*,\e)$ where $c_*$ is 
the multiset consisting of the sizes of the Jordan blocks of an element of $\g$; $\e(j)$ is equal to $0$ if 
$((g-1)^{j-1}(x),x)=0$ for all $x\in\ker(g-1)^j$ ($g\in\g$) and $\e(j)=1$ otherwise (see \cite{\SPA}).

If $Q\ne0,\k=1,p\ne2$ we identify $\uuG=\cq_\nn$ by associating to $\g\in\uuG$ the multiset consisting of the 
sizes of the Jordan blocks of an element of $\g$. 

If $Q\ne0,\k=1,p=2$ we identify $\uuG=\ct_{2n}^{(2)}$ by associating to $\g\in\uuG$ the pair $(c_*,\e)$
corresponding as above to the image of $\g$ under the obvious bijective homomorphism from $G$ to a symplectic 
group of an $\nn-1$ dimensional vector space.

If $Q\ne0,\k=0$ we denote by $\uuG_{(0)}$ (resp. $\uuG_{(1)}$) the set of unipotent classes in $G$ which are not
conjugacy classes in $Is(V)$ (resp. are also conjugacy classes in $Is(V)$); let $\ti{\uuG}$ be the set of
$Is(V)$-conjugacy classes of unipotent elements of $G$. We have an obvious imbedding $\uuG_{(1)}\sub\ti{\uuG}$.

If $Q\ne0,\k=0,p\ne2$ we identify $\ti{\uuG}=\cq_\nn$ by associating to $\g\in\ti{\uuG}$ the multiset 
consisting of the sizes of the Jordan blocks of an element of $\g$. 

If $Q\ne0,\k=0,p=2$ we identify $\ti{\uuG}=\ti{\ct}_{2n}^{(2)}$  by associating to $\g\in\ti{\uuG}$ the pair
$(c_*,\e)\in\ct_{2n}^{(2)}$ corresponding as above to the image of $\g$ under the obvious imbedding of $Is(V)$ 
into the symplectic group of $V,(,)$.

We define $\io:\ca^1_{2n}@>>>\ct_{2n}$ by $(r_*,p_*)\m c_*$ where the multiset of entries of
$c_*$ is the union of the multiset of entries of $r_*$ with the multiset of entries of $p_*$.
We define $\io^{(2)}:\ca^1_{2n}@>>>\ct_{2n}^{(2)}$ by $(r_*,p_*)\m(c_*,\e)$ where 
$c_*=\io(r_*,p_*)$ and for any $j\in\{2,4,6,\do\}$ such that $\mu_j(c_*)\in\{2,4,6,\do\}$, we have
$\e(j)=1$ if $j=r_i$ for some $i$ and $\e(j)=0$, otherwise.
When $\k=0$ we define $\ti{\io}^{(2)}:\ca^0_{2n}@>>>\ti{\ct}_{2n}^{(2)}$ to be the restriction of 
$\io^{(2)}:\ca^1_{2n}@>>>\ct_{2n}^{(2)}$.

In the remainder of this section we assume that $\nn\ge3$. We define $\X:\cs^\k@>>>\cp^1$ by
$$(r_1\ge r_2\ge\do\ge r_\s)\m(r_1+\ps(1)\ge r_2+\ps(2)\ge\do\ge r_\s+\ps(\s))$$
if $\s+\k$ is even,
$$(r_1\ge r_2\ge\do\ge r_\s)\m(r_1+\ps(1)\ge r_2+\ps(2)\ge\do\ge r_\s+\ps(\s)\ge1)$$
if $\s+\k$ is odd, where $\ps:[1,\s]@>>>\{-1,0,1\}$ is as follows:

if $t\in[1,\s]$ is odd and $r_t<r_x$ for any $x\in[1,t-1]$ then $\ps(t)=1$;

if $t\in[1,\s]$ is even and $r_x<r_t$ for any $x\in[t+1,\s]$, then $\ps(t)=-1$;

for all other $t\in[1,\s]$ we have $\ps(t)=0$.

We define $\io':\ca^\k_{2n}@>>>\cq_\nn$ by $(r_*,p_*)\m c_*$ where the multiset of entries of is the union of the 
multiset of entries of $\X(r_*)$ with the multiset of entries of $p_*$. (We will see below that $\io'$ is well 
defined.)

With the identifications above and those in 1.3, we see that the map $\Ph:\uWW@>>>\uuG$ becomes:

(a) $\io:\ca^1_{2n}@>>>\ct_{2n}$ if $Q=0,p\ne2$, see \cite{\WE, 3.7, 1.1};

(b) $\io':\ca^1_{2n}@>>>\cq_\nn$   if $Q\ne0,\k=1,p\ne2$, see \cite{\WE, 3.8, 1.1};

(c) $\io^{(2)}:\ca^1_{2n}@>>>\ct_{2n}^{(2)}$ if $Q=0,p=2$, see \cite{\WE, 4.6, 1.1};

(d) $\io^{(2)}:\ca^1_{2n}@>>>\ct_{2n}^{(2)}$ if $Q\ne0,\k=1,p=2$, see \cite{\WE, 4.6, 1.1};
\nl
and that (when $Q\ne0,\k=0$) the map $\ti\Ph:\ti{\uWW}@>>>\ti{\uuG}$ induced by $\Ph$ becomes:

(e) $\io':\ca^0_{2n}@>>>\cq_\nn$ if $Q\ne0,\k=0,p\ne2$, see \cite{\WE, 3.9, 1.1};

(f) $\ti{\io}^{(2)}:\ca^0_{2n}@>>>\ti{\ct}_{2n}^{(2)}$ if $Q\ne0,\k=0,p=2$, see \cite{\WE, 4.6, 1.1}.
\nl
In particular, $\io'$ is well defined.

We see that to prove 0.2 it is enough to prove the following statement:

(g) {\it In each of the cases (a)-(f), the function $(r_*,p_*)\m\t_{p_*}/2$ on any fibre of the map described in 
that case, reaches its minimum value at exactly one element in that fibre.}
\nl
We have used that in the case where $Q\ne0,\k=0$, the fibre of $\Ph$ over any element in $\uuG_{(0)}$ has 
exactly one element (necessarily in $\uWW_0$) and the fibre of $\Ph$ over any element in $\uuG_{(1)}$ is 
contained in $\uWW_1$ and is the same as the fibre of $\ti\Ph$ over that element.

\subhead 1.5\endsubhead
We prove 1.4(g) in the case 1.4(a). Let $c_*\in\ct_{2n}$. Let $(r_*,p_*)\in\io\i(c_*)$. Let $M_e=\mu_e(r_*)$,
$N_e=\mu_e(p_*)$, $Q_e=\mu_e(c_*)$ so that $M_e+N_e=Q_e$. If $e$ is odd then $M_e=0$ hence $N_e=Q_e$. Thus 
$\sum_eN_e\ge\sum_{e\text{ odd}}Q_e$. We see that the minimum value of the function 1.4(g) on $\io\i(c_*)$ is 
reached when $M_e=Q_e,N_e=0$ for $e$ even and $M_e=0,N_e=Q_e$ for $e$ odd. This proves 1.4(g) in our case.

\subhead 1.6\endsubhead
We prove 1.4(g) in the case 1.4(c). Let $(c_*,\e)\in\ct_{2n}^{(2)}$. Let $(r_*,p_*)\in(\io^{(2)})\i(c_*,\e)$. Let 
$M_e=\mu_e(r_*)$, $N_e=\mu_e(p_*)$, $Q_e=\mu_e(c_*)$ so that $M_e+N_e=Q_e$. If $e$ is odd then $M_e=0$ hence 
$N_e=Q_e$. If $e$ is even, $Q_e$ is even $\ge2$ and $\e(e)=0$ then $M_e=0$. Thus
$$\sum_eN_e\ge\sum_{e\text{ odd}}Q_e+\sum_{e\text{ even},Q_e\text{ even},\ge2,\e(e)=0}Q_e.$$
We see that the minimum value of the function 1.4(g) on $(\io^{(2)})\i(c_*,\e))$ is reached when $M_e=0,N_e=Q_e$ 
for $e$ odd and for $e$ even with $Q_e\text{ even},\ge2,\e(e)=0$ and $M_e=Q_e,N_e=0$ for all other $e$. This 
proves 1.4(g) in our case. The same proof yields 1.4(g) in the case 1.4(d); a similar proof yields 1.4(g) in the 
case 1.4(f).

\subhead 1.7\endsubhead
In this subsection we prove 1.4(g) in the cases 1.4(b) and 1.4(e).

Let $\car$ be the set of all $r_*=(r_1\ge r_2\ge\do\ge r_\t)\in\cq$ such that the following conditions are
satisfied. Let $J_{r_*}=\{k\in[1,\t];r_k\text{ is odd}\}$. We write the multiset $\{r_k;k\in J_{r_*}\}$ as a 
sequence $r^1\ge r^2\ge\do\ge r^s$. (We have necessarily $\t=s\mod2$.) We require that:

-if $\t\ne0$ then $1\in J_{r_*}$;

-if $\t\ne0$ is even then $\t\in J_{r_*}$;

-if $u\in[1,s-1]$ is odd then $r^u>r^{u+1}$;

-if $u\in[1,s-1]$ is even then there is no $k'\in[1,\t]$ such that $r^u>r_{k'}>r^{u+1}$.
\nl
For $N\in\NN$ we set $\car_N=\car\cap\cq_N$.

Note that for any even $N\in\NN$, $\X$ in 1.1 restricts to a bijection $\cs_N^\k@>>>\car_{N+\k}$ 
with inverse map $\car_{N+\k}@>>>\cs_N^\k$ given by
$$(r_1\ge r_2\ge\do\ge r_\t)\m(r_1+\z(1)\ge r_2+\z(2)\ge\do\ge r_\t+\z(\t))$$
if $r_\t>\k$,
$$(r_1\ge r_2\ge\do\ge r_\t)\m(r_1+\z(1)\ge r_2+\z(2)\ge\do\ge r_{\t-1}+\z(\t-1))$$
if $r_\t=\k$, where $\z:[1,\t]@>>>\{-1,0,1\}$ is given by $\z(k)=(-1)^k(1-(-1)^{r_k})/2$. (Thus, 
$\z(k)=(-1)^k$ if $r_k$ is odd and $\z(k)=0$ if $r_k$ is even. We have $r_k+\z(k)\in2\NN$ for any $k$ and 
$r_k+\z(k)\ge r_{k+1}+\z(k+1)$ for $k\in[1,\t-1]$.)

Let $A_\nn$ be the set of all pairs $(r_*,p_*)\in\car\T\tcp$ such that $|r_*|+|p_*|=\nn$. We define 
$\ti\io:A_\nn@>>>\cq_\nn$ by $(r_*,p_*)\m c_*$ where the multiset of entries of $c_*$ is the union of the 
multiset of entries of $r_*$ with the multiset of entries of $p_*$. In view of the bijection  
$\cs_N^\k@>>>\car_{N+\k}$ defined by the restriction of $\Xi$ we see that to prove 1.4(g) in our case it is 
enough to prove the following statement.

(a) {\it For any $c_*\in\cq_\nn$ there is exactly one element $(r_*,p_*)\in\ti\io\i(c_*)$ such that the number of
entries of $p_*$ is minimal.}
\nl
Let $c_*=(c_1\ge c_2\ge\do\ge c_\t)\in\cq_\nn$. Let $K=\{k\in[1,\t];c_k\text{ is odd}\}$. We write the multiset 
$\{c_k;k\in K\}$ as a sequence $c^1\ge c^2\ge\do\ge c^t$. (We have necessarily $\t=\nn=t\mod2$.) We associate to 
$c_*$ an element $(r_*,p_*)\in\cq\T\ti\cp$ by specifying $M_e=\mu_e(r_*)$, $N_e=\mu_e(p_*)$ for $e\ge1$. Let 
$Q_e=\mu_e(c_*)$. 

(i) If $e\in2\NN+1$ and $Q_e=2g+1$, then $M_e=1$, $N_e=2g$.

(ii) If $e\in2\NN+1$ and $Q_e=2g$, so that $c^d=c^{d+1}=\do=c^{d+2g-1}=e$ with $d$ even, then $M_e=2$, $N_e=2g-2$
(if $g>0$) and $M_e=N_e=0$ (if $g=0$).

(iii) If $e\in2\NN+1$ and $Q_e=2g$ so that $c^d=c^{d+1}=\do=c^{d+2g-1}=e$ with $d$ odd, then $M_e=0$, $N_e=2g$.

Thus the odd entries of $r_*$ are defined. We write them in a sequence $r^1\ge r^2\ge\do\ge r^s$.

(iv) If $e\in2\NN+2$, $Q_e=2g$ and if 

($*$) $r^{2v}>e>r^{2v+1}$ for some $v$, or $e>r^1$, or $r^s>e$ (with $s$ even),
\nl
then $M_e=0,N_e=2g$. 

(v) If $e\in2\NN+2$, $Q_e=2g$  and if ($*$) does not hold, then $M_e=2g,N_e=0$. 

Now $r_*\in\cq,p_*\in\tcp$ are defined and $|r_*|+|p_*|=\nn$.

Assume that $|r_*|>0$; then from (iv) we see that the largest entry of $r_*$ is odd. Assume that $|r_*|>0$ and
$\nn$ is even; then from (iv) we see that the smallest entry of $r_*$ is odd. If $u\in[1,s-1]$ and $r^u=r^{u+1}$ 
then from (i),(ii),(iii) we see that $u$ is even. If $u\in[1,s-1]$ and there is $k'\in[1,\t]$ such that 
$r^u>r_{k'}>r^{u+1}$, then $r_{k'}$ is even and $e=r_{k'}$ is as in (v) and $u$ must be odd. We see that 
$r_*\in\car$.

We see that $c_*\m(r_*,p_*)$ is a well defined map $\ti\io':\cq_\nn@>>>A_\nn$; moreover, 
$\ti\io\ti\io':\cq_\nn@>>>\cq_\nn$ is the identity map. 

We preserve the notation for $c_*,r_*,p_*$ as above (so that $(r_*,p_*)\in\ti\io\i(c_*)$) and we assume that
$(r'_*,p'_*)\in\ti\io\i(c_*)$. We write the odd entries of $r'_*$ in a sequence 
$r'{}^1\ge r'{}^2\ge\do\ge r'{}^{s'}$.

Let $M'_e=\mu_e(r'_*)$, $N'_e=\mu_e(p'_*)$ for $e\ge1$. Note that $M'_e+N'_e=M_e+N_e$.

In the setup of (i) we have $M'_e=1,N'_e=N_e$. (Indeed, $M'_e+N'_e$ is odd, $N'_e$ is even hence $M'_e$ is odd. 
Since $M'_e$ is $0,1$ or $2$ we see that it is $1$.)

In the setup of (ii) and assuming that $g>0$ we have $M'_e=2,N'_e=N_e$ or $M'_e=0,N'_e=N_e+2$. (Indeed, 
$M'_e+N'_e$ is even, $N'_e$ is even hence $M'_e$ is even. Since $M'_e$ is $0,1$ or $2$ we see that it is $0$ or 
$2$.) If $g=0$ we have $M'_e=N'_e=0$.

In the setup of (iii) we have $M'_e=0,N'_e=N_e$. (Indeed, $M'_e+N'_e$ is even, $N'_e$ is even hence $M'_e$ is
even. Since $M'_e$ is $0,1$ or $2$ we see that it is $0$ or $2$. Assume that $M'_e=2$. Then $e=r'{}^u=r'{}^{u+1}$
with $u$ even in $[1,s'-1]$. We have $c^d=c^{d+1}=\do=c^{d+2g-1}=e$ with $d$ odd. From the definitions we see
that $u=d\mod2$ and we have a contradiction. Thus, $M'_e=0$.)

Now the sequence $r'{}^1\ge r'{}^2\ge\do\ge r'{}^{s'}$ is obtained from the sequence $r^1\ge r^2\ge\do\ge r^s$ by
deleting some pairs of the form $r^{2h}=r^{2h+1}$. Hence in the setup of (iv) we have 
$r^{2v}=r'{}^{2v'}>e>r'{}^{2v'+1}=r^{2v+1}$ for some $v'$ or $e>r'{}^1$ or $r'{}^{s'}>e$ (with $s,s'$ even) and 
we see that $M'_e=0$ so that $N'_e=2g=N_e$. 

In the setup of (v) we have $N'_e\ge 0$.

We see that in all cases we have $N'_e\ge N_e$. It follows that $\sum_eN'_e\ge\sum_eN_e$ (and the equality 
implies that $N'_e=N_e$ for all $e$ hence $(r'_*,p'_*)=(r_*,p_*)$). This proves (a) and completes the proof of 
1.4(g) in all cases hence the proof of 0.2 for $G$ almost simple of type $B,C$ or $D$.

\head 2. Exceptional groups\endhead
\subhead 2.1\endsubhead
In 2.2-2.6 we describe explicitly the map $\Ph:\uWW@>>>\uuG$ in the case where $G$ is a simple exceptional group 
in the form of tables. Each table consists of lines of the form $[a,b,\do,r]\m\g$ where $\g\in\uuG$ is specified 
by its name in \cite{\SPAE} and $a,b,\do,r$ are the elements of $\uWW$ which are mapped by $\Ph$ to $\g$ (here 
$a,b,\do,r$ are specified by their name in \cite{\CA}); by inspection we see that 0.2 holds in each case and in 
fact $\Ps(\g)=a$ is the first element of $\uWW$ in the list $a,b,\do,r$. The tables are obtained from the results
in \cite{\WE}.

\subhead 2.2. Type $G_2$\endsubhead
If $p\ne3$ we have

$[A_0]\m A_0$

$[A_1]\m A_1$ 

$[A_1+\tA_1,\tA_1] \m \tA_1$  

$[A_2]\m G_2(a_1)$

$[G_2]\m G_2$.
\nl
When $p=3$ the line $[A_1+\tA_1,\tA_1]\m\tA_1$ should be replaced by $[A_1+\tA_1]\m\tA_1$, $[\tA_1]\m(\tA_1)_3$. 

\subhead 2.3. Type $F_4$\endsubhead
If $p\ne2$ we have 

$[A_0]\m A_0$

$[A_1]\m A_1$

$[2A_1,\tA_1]\m \tA_1$   

$[4A_1,3A_1,2A_1+\tA_1,A_1+\tA_1]\m A_1+\tA_1$ 

$[A_2]\m A_2$

$[\tA_2]\m\tA_2$

$[A_2+\tA_1]\m A_2+\tA_1$

$[A_2+\tA_2,\tA_2+A_1]\m\tA_2+A_1$ 

$[A_3,B_2]\m B_2$  

$[A_3+\tA_1,B_2+A_1]\m C_3(a_1)$  

$[D_4(a_1)]\m F_4(a_3)$

$[D_4,B_3]\m B_3$ 

$[C_3+A_1,C_3]\m C_3$

$[F_4(a_1)]\m F_4(a_2)$ 

$[B_4]\m F_4(a_1)$

$[F_4]\m F_4$.
\nl
When $p=2$ the lines $[2A_1,\tA_1]\m\tA_1$, $[A_2+\tA_2,\tA_2+A_1]\m\tA_2+A_1$, $[A_3,B_2]\m B_2$,
$[A_3+\tA_1,B_2+A_1]\m C_3(a_1)$, should be replaced by

$[2A_1]\m\tA_1$, $[\tA_1]\m(\tA_1)_2$

$[A_2+\tA_2]\m\tA_2+A_1$, $[\tA_2+A_1]\m(\tA_2+A_1)_2$

$[A_3]\m B_2$, $[B_2]\m(B_2)_2$

$[A_3+\tA_1]\m C_3(a_1)$, $[B_2+A_1]\m(C_3(a_1))_2$
\nl
respectively.

\subhead 2.4. Type $E_6$\endsubhead
We have

$[A_0]\m A_0$ 

$[A_1]\m A_1$

$[2A_1]\m2A_1$ 

$[4A_1,3A_1]\m3A_1$

$[A_2]\m A_2$

$[A_2+A_1]\m A_2+A_1$ 

$[2A_2]\m2A_2$

$[A_2+2A_1]\m A_2+2A_1$

$[A_3]\m A_3$

$[3A_2,2A_2+A_1]\m 2A_2+A_1$ 

$[A_3+2A_1,A_3+A_1]\m A_3+A_1$ 

$[D_4(a_1)]\m D_4(a_1)$ 

$[A_4]\m A_4$ 

$[D_4]\m D_4$

$[A_4+A_1]\m A_4+A_1$    

$[A_5+A_1,A_5]\m A_5$

$[D_5(a_1)]\m D_5(a_1)$ 

$[E_6(a_2)]\m A_5+A_1$ 

$[D_5]\m D_5$ 

$[E_6(a_1)]\m E_6(a_1)$ 

$[E_6]\m E_6$.

\subhead 2.5. Type $E_7$\endsubhead
If $p\ne 2$ we have  

$[A_0]\m A_0$

$[A_1]\m A_1$ 

$[2A_1]\m 2A_1$

$[(3A_1)']\m(3A_1)''$ 

$[(4A_1)'',(3A_1)'']\m(3A_1)'$

$[A_2]\m A_2$  

$[7A_1,6A_1,5A_1,(4A_1)']\m 4A_1$ 

$[A_2+A_1]\m A_2+A_1$ 

$[A_2+2A_1]\m A_2+2A_1$ 

$[A_3]\m A_3$ 

$[2A_2]\m 2A_2$
 
$[A_2+3A_1]\m A_2+3A_1$  

$[(A_3+A_1)']\m(A_3+A_1)''$ 

$[3A_2,2A_2+A_1]\m 2A_2+A_1$ 

$[(A_3+2A_1)'',(A_3+A_1)'']\m(A_3+A_1)'$ 

$[D_4(a_1)]\m D_4(a_1)$ 

$[A_3+3A_1,(A_3+2A_1)']\m A_3+2A_1$ 

$[D_4]\m D_4$ 

$[D_4(a_1)+A_1]\m D_4(a_1)+A_1$ 

$[D_4(a_1)+2A_1,A_3+A_2]\m A_3+A_2$ 

$[2A_3+A_1,A_3+A_2+A_1]\m A_3+A_2+A_1$ 

$[A_4]\m A_4$ 

$[D_4+3A_1,D_4+2A_1,D_4+A_1]\m D_4+A_1$ 

$[A'_5]\m A''_5$ 

$[A_4+A_1]\m A_4+A_1$ 

$[D_5(a_1)]\m D_5(a_1)$  

$[A_4+A_2]\m A_4+A_2 $

$[(A_5+A_1)'',A''_5]\m A'_5$  

$[A_5+A_2,(A_5+A_1)']\m(A_5+A_1)''$  

$[D_5(a_1)+A_1]\m D_5(a_1)+A_1$ 

$[E_6(a_2)]\m(A_5+A_1)'$ 

$[D_6(a_2)+A_1,D_6(a_2)]\m D_6(a_2) $ 

$[E_7(a_4)]\m D_6(a_2)+A_1$ 

$[D_5]\m D_5$   

$[A_6]\m A_6$

$[D_5+A_1]\m D_5+A_1$ 

$[D_6(a_1)]\m D_6(a_1)$ 

$[A_7]\m D_6(a_1)+A_1$ 

$[E_6(a_1)]\m E_6(a_1)$  

$[D_6+A_1,D_6]\m D_6$ 

$[E_6]\m E_6$ 

$[E_7(a_3)]\m D_6+A_1$  

$[E_7(a_2)]\m E_7(a_2)$ 

$[E_7(a_1)]\m E_7(a_1)$ 

$[E_7]\m E_7$.
\nl
If $p=2$, the line $[D_4(a_1)+2A_1,A_3+A_2]\m A_3+A_2$ should be replaced by $[D_4(a_1)+2A_1]\m A_3+A_2$, 
$[A_3+A_2]\m (A_3+A_2)_2$.

\subhead 2.6. Type $E_8$ \endsubhead
If $p\ne2,3$ we have

$[A_0]\m A_0$ 

$[A_1]\m A_1$ 

$[2A_1]\m 2A_1$

$[(4A_1)',3A_1]\m 3A_1 $

$[A_2]\m A_2$

$[8A_1,7A_1,6A_1,5A_1,(4A_1)'']\m 4A_1$

$[A_2+A_1]\m A_2+A_1$

$[A_2+2A_1]\m A_2+2A_1$

$[A_3]\m A_3$

$[A_2+4A_1,A_2+3A_1]\m A_2+3A_1$

$[2A_2]\m 2A_2$

$[3A_2,2A_2+A_1]\m 2A_2+A_1$

$[(A_3+2A_1)',A_3+A_1]\m A_3+A_1$

$[D_4(a_1)]\m D_4(a_1)$

$[4A_2,3A_2+A_1,2A_2+2A_1]\m 2A_2+2A_1$

$[D_4]\m D_4$

$[A_3+4A_1,A_3+3A_1,(A_3+2A_1)'']\m A_3+2A_1$

$[D_4(a_1)+A_1]\m D_4(a_1)+A_1$

$[(2A_3)',A_3+A_2\m A_3+A_2$  

$[A_4]\m A_4$

$[2A_3+2A_1,A_3+A_2+2A_1,2A_3+A_1,A_3+A_2+A_1]\m A_3+A_2+A_1$

$[D_4(a_1)+A_2]\m D_4(a_1)+A_2$

$[D_4+4A_1,D_4+3A_1,D_4+2A_1,D_4+A_1]\m D_4+A_1$

$[2D_4(a_1),D_4(a_1)+A_3,(2A_3)'']\m 2A_3$

$[A_4+A_1]\m A_4+A_1$

$[D_5(a_1)]\m D_5(a_1)$

$[A_4+2A_1]\m A_4+2A_1$

$[A_4+A_2]\m A_4+A_2$

$[A_4+A_2+A_1]\m A_4+A_2+A_1$

$[D_5(a_1)+A_1]\m D_5(a_1)+A_1$

$[(A_5+A_1)',A_5]\m A_5$

$[D_4+A_3,D_4+A_2]\m D_4+A_2$ 

$[E_6(a_2)]\m (A_5+A_1)''$

$[2A_4,A_4+A_3]\m A_4+A_3$

$[D_5]\m D_5$

$[D_5(a_1)+A_3,D_5(a_1)+A_2]\m D_5(a_1)+A_2$

$[A_5+A_2+A_1,A_5+A_2,A_5+2A_1,(A_5+A_1)'']\m (A_5+A_1)'$

$[E_6(a_2)+A_2,E_6(a_2)+A_1]\m A_5+2A_1$

$[2D_4,D_6(a_2)+A_1,D_6(a_2)]\m D_6(a_2)$

$[E_7(a_4)+A_1,E_7(a_4)]\m A_5+A_2$      

$[D_5+2A_1,D_5+A_1]\m D_5+A_1$

$[E_8(a_8)]\m2A_4$

$[D_6(a_1)\m D_6(a_1)]$

$[A_6]\m A_6$

$[A_6+A_1]\m A_6+A_1$

$[A'_7]\m D_6(a_1)+A_1$

$[A_7+A_1,D_5+A_2]\m D_5+A_2$ 

$[E_6(a_1)]\m E_6(a_1)$

$[D_6+2A_1,D_6+A_1,D_6]\m D_6$

$[D_7(a_2)]\m D_7(a_2)$

$[E_6]\m E_6$

$[D_8(a_3),A''_7]\m A_7$                      

$[E_6(a_1)+A_1]\m E_6(a_1)+A_1 $

$[E_7(a_3)]\m D_6+A_1$

$[A_8]\m D_8(a_3)$

$[D_8(a_2),D_7(a_1)]\m D_7(a_1)$ 

$[E_6+A_2,E_6+A_1]\m E_6+A_1$

$[E_7(a_2)+A_1,E_7(a_2)]\m E_7(a_2)$

$[E_8(a_6)]\m A_8$

$[E_8(a_7)]\m E_7(a_2)+A_1$

$[D_8(a_1),D_7]\m D_7$

$[E_8(a_3)]\m D_8(a_1)$

$[E_7(a_1)]\m E_7(a_1)$

$[D_8]\m E_7(a_1)+A_1$

$[E_8(a_5)]\m D_8$

$[E_7+A_1,E_7]\m E_7$

$[E_8(a_4)]\m E_7+A_1$

$[E_8(a_2)]\m E_8(a_2)$

$[E_8(a_1)]\m E_8(a_1)$

$[E_8]\m E_8$.
\nl
If $p=3$ the line $[D_8(a_3),A''_7]\m A_7$ should be replaced by $[D_8(a_3)]\m A_7$, $[A''_7]\m(A_7)_3$. If $p=2$
the lines $[(2A_3)',A_3+A_2]\m A_3+A_2$,  $[D_4+A_3,D_4+A_2]\m D_4+A_2$, 
$[A_7+A_1,D_5+A_2]\m D_5+A_2$, $[D_8(a_2),D_7(a_1)]\m D_7(a_1)$ should be replaced by    

$[(2A_3)']\m A_3+A_2$, $[A_3+A_2]\m(A_3+A_2)_2$

$[D_4+A_3]\m D_4+A_2$,  $[D_4+A_2]\m (D_4+A_2)_2$
 
$[A_7+A_1]\m D_5+A_2$, $[D_5+A_2]\m (D_5+A_2)_2$
   
$[D_8(a_2)]\m D_7(a_1)$, $[D_7(a_1)]\m (D_7(a_1))_2$    
\nl
respectively.

\subhead 2.7\endsubhead
From the tables above we see that (assuming that $G$ is almost simple of exceptional type) the following holds.

(a) {\it If $\g$ is a distinguished unipotent class in $G$ then $\Ph\i(\g)$ consists of a single conjugacy class 
in $\WW$.}
\nl
In fact (a) is also valid without any assumption on $G$. Indeed, assume that $C\in\Ph\i(\g)$. From the arguments 
in \cite{\WE, 1.1} we see that if $C$ is not elliptic then $\Ph(C)$ is not distinguished. Thus $C$ must be 
elliptic. Then the desired result follows from the injectivity of the restriction of $\Ph$ to elliptic conjugacy 
classes in $\WW$, see \cite{\WE, 0.6}.

\subhead 2.8\endsubhead
We have the following result (for general $G$).

(a) {\it If $C$ is an elliptic conjugacy class in $\WW$, then $C=\Ps(\Ph(C))$. In particular, $C$ is in the image
of $\Ps:\uuG@>>>\uWW$.}
\nl
Since $C\in\Ph\i(\Ph(C))$, (a) follows immediately from Theorem 0.2.

\head 3. Proof of Theorem 0.4\endhead
\subhead 3.1\endsubhead
There is a well defined (injective) map $\p':\uuG_0@>>>\uuG$,  $\g_0\m\p'(\g_0)$, where $\p'(\g_0)$ is the unique 
unipotent class in $G$ which has the same Springer representation of $\WW$ as $\g_0$. One can show that $\p'$ 
coincides with the order preserving and dimension preserving imbedding defined in \cite{\SPA, III,5.2}. 

To prove Theorem 0.4(a) we can assume that $G,G_0$ are almost simple. It is also enough to prove the theorem for
a single $G$ in each isogeny class of almost simple groups. We can assume that $p$ is a bad prime for $G$ (if 
$p$ is not a bad prime, the result is obvious). Now $G,G_0$ cannot be of type $A$ since $p$ is a bad prime 
for $G$. In the case where $G$ is of exceptional type, the theorem follows by inspection of the tables in Section
2. In the case where $G$ is of type $B,C$ or $D$ so that $p=2$, we define $\ti\r:\uuG@>>>\uuG_0$ by 
$\ti\r(\g)=\g_0$ where $\g\in\uuG$ is in the "unipotent piece" of $G$ indexed by $\g_0\in\uuG_0$ (in the sense of
\cite{\USCIII}). It is enough to prove that
$$\ti\r\Ph=\Ph_0.\tag a$$
(This would prove the existence of $\r$ in 0.4 and that $\r=\ti\r$.) If $G,G_0$ are of type 
$C_n$ ($n\ge2$) then 
$\Ph_0,\Ph$ may be identified with $\io:\ca^1_{2n}@>>>\ct_{2n}$, $\io^{(2)}:\ca^1_{2n}@>>>\ct_{2n}^{(2)}$ (see 
1.4) and by \cite{\USCI}, $\ti\r$ may be identified with $\ct_{2n}^{(2)}@>>>\ct_{2n}$, $(c_*,\e)\m c_*$ (notation
of 1.4). Then the identity (a) is obvious. The proof of (a) for $G$ of type $B$ and $D$ is given in 3.5.

Similarly to prove Theorem 1.4(b) it is enough to prove that
$$\Ps\p'=\Ps_0.\tag b$$
(This would prove the existence of $\p$ in 0.4 and that $\p=\p'$.) Again it is enough to prove (b) in the case 
where $G$ is of type $B,C$ or $D$ so that $p=2$. The proof is given in 3.9.

\subhead 3.2\endsubhead
In this subsection we assume that all simple factors of $G$ are of type $A,B,C$ or $D$. Then
$\ti\r:\uuG@>>>\uuG_0$ is defined as in 3.1. We show:

(a) {\it $\Ph_0(C)=\ti\r\Ph(C)$ for any elliptic conjugacy class $C$ in $\WW$.}
\nl
From the definitions we see that

(b) $\ti\r\p'=1$.
\nl
Now for $C$ as in (a) we have by definition $\Ph(C)=\p'\Ph_0(C)$. To prove (a) it is enough to show that 
$\Ph_0(C)=\ti\r\p\Ph_0(C)$. But this follows from (b).

\subhead 3.3\endsubhead
In this subsection we assume that $V,Q,(,),\nn,Is(V)$ are as in 1.2 with $p=2$ and that either $Q\ne0$ or $V=0$. 
Let $SO(V)$ be the identity component of $Is(V)$. Let $U_V$ be the set of unipotent elements in $SO(V)$. 

We say that $u\in U_V$ is split if the corresponding pair $(c_*,\e)$ (see 1.4) satisfies $\mu_j(c_*)=$even for 
all $j$ and $\e(j)=0$ for all even $j$ such that $\mu_j(c_*)>0$ ($\mu_j$ as in 1.1).

For $u\in U_V$ let $e=e_u$ be the smallest integer $\ge0$ such that $(u-1)^e=0$. (When $u=1$ we have $e=1$ if 
$\nn>0$ and $e=0$ if $\nn=0$. When $u\ne1$ we have $e\ge2$.) If $u\ne1$ we define $\l=\l_u:V@>>>\kk$ by 
$\l(x)=\sqrt{(x,(u-1)^{e-1}x)}$, a linear form on $V$; in this case we define $L=L_u$ as follows:
$$\align&L=(u-1)^{e-1}V\text{ if }\l=0;\qua L=(\ker\l)^\pe\text{ if }\l\ne0,\nn=\text{even};\\&
L=\{x\in(\ker\l)^\pe;Q(x)=0\text{ if }\l\ne0,\nn=\text{odd}.\endalign$$
Note that $L\ne0$.

As in \cite{\USCII, 2.5}, for any $u\in U_V$ we define subspaces $V^a=V^a_u$ ($a\in\ZZ$) of $V$ as follows. If 
$u=1$ we set $V^a=V$ for $a\le0$, $V^a=0$ for $a\ge1$. If $u\ne1$ so that $e=e_u\ge2$, and $\l=\l_u$, $L=L_u$ are
defined, we have $L\sub L^\pe$, $Q|L=0$ (see \cite{\USCII,2.2(a),(e)}) and we set $\bV=L^\pe/L$; now $\bV$ has an
induced nondegenerate quadratic form and there is an induced unipotent element $\bar u\in SO(\bV)$; let
$r:L^\pe@>>>L^\pe/L=\bV$ be the canonical map. Since $\dim\bV<\dim V$, we can assume by induction that
$\bV^a=\bV^a_{\bar u}$ are defined for $a\in\ZZ$. We set 
$$V^a=V\text{ if }a\le1-e,\qua V^a=r\i(\bV^a)\text{ if }2-e\le a\le e-1,\qua V^a=0\text{ if }a\ge e$$
when $\l=0$ and
$$V^a=V\text{ if }a\le-e,\qua V^a=r\i(\bV^a)\text{ if }1-e\le a\le e,\qua V^a=0\text{ if }a\ge1=e$$
when $\l\ne0$. This completes the inductive definition of $V^a=V^a_u$.

Note that if $u\in U_V$ is split and $u\ne1$ then $\bar u$ (as above) is split.

Now assume that $V=V'\op V''$ where $V',V''$ are subspaces of $V$ such that $(V',V'')=0$ so that $Q|_{V'}$, 
$Q|_{V'}$ are nondegenerate and let $u\in U_V$ be such that $V',V''$ are $u$-stable; we set $u'=u|_{V'}$, 
$u''=u|_{V''}$. We assume that $u'\in U_{V'}$, $u''\in U_{V''}$. For $a\in\ZZ$ we set $V^a=V^a_u$, 
$V'{}^a=V'{}^a_{u'}$, $V''{}^a=V''{}^a_{u''}$. We show:

(a) {\it if $u''$ is split, then $V^a=V'{}^a\op V''{}^a$ for all $a$.}
\nl
If $u=1$, then (a) is trivial. So we can assume that $u\ne1$. We can also assume that (a) is true when $V$ is 
replaced by a vector space of smaller dimension. Let $e=e_u,e'=e_{u'},e''=e_{u''}$. Let $\l=\l_u,L=L_u$. Let
$r:L^\pe@>>>L^\pe/L=\bV$ be as above. If $u'\ne1$ (resp. $u''\ne1$) we define $\l',L',r',\bV'$ (resp. 
$\l'',L'',r'',\bV''$) in terms of $u',Q_{V'}$ (resp. $u'',Q_{V''}$) in the same way as $\l,L,r,\bV$ were defined 
in terms of $u,Q$. Note that $\l''=0$ (when $u''\ne1$) since $u''$ is split.

Assume first that $e''>e'$. We have $e=e''$ hence $u''\ne1$. Moreover, $L=0\op L''$, $\bV=V'\op\bV''$. By the
induction hypothesis we have $\bV^a=V'{}^a\op\bV''{}^a$ for all $a$. If $a\le1-e$ then 
$V^a=V=V'\op V''=V'{}^a\op V''{}^a$. (We use that $V'=V'{}^a$; if $e'\ge2$ this follows from $a\le-e'$; if 
$e'\le1$ this follows from $a\le0$). If $2-e\le a\le e-1$ we have 
$$V^a=(0\op r'')\i(\bV^a)=(0\op r'')\i(V'{}^a\op\bV''{}^a)=V'{}^a\op V''{}^a$$
(we use that $e''=e$); if $a\ge e$ then $V^a=0=V'{}^a\op V''{}^a$. (We use that $V'{}^a=0$; if $e'\ge2$ this
follows from $a\ge e'+1$; if $e'\le1$ this follows from $a\ge1$.)

Next we assume that $e'=e''$ (hence both are equal to $e\ge2$) and that $\l=0$ (hence $\l'=0$). Then 
$L=L'\op L''$, $\bV=\bV'\op\bV''$. By the induction hypothesis we have $\bV^a=\bV'{}^a\op\bV''{}^a$ for all $a$.
If $a\le1-e$ then $V^a=V=V'\op V''=V'{}^a\op V''{}^a$. If $2-e\le a\le e-1$, we have 
$$V^a=(r'\op r'')\i)(\bV^a)=(r'\op r'')\i)(\bV'{}^a\op\bV''{}^a)=V'{}^a\op V''{}^a.$$
If $a\ge e$ then $V^a=0=V'{}^a\op V''{}^a$. 

Next we assume that $e'>e''$ (hence $e'=e\ge2$) and $\l=0$ (hence $\l'=0$). Then $L=L'\op 0$, $\bV=\bV'\op V''$.
By the induction hypothesis we have $\bV^a=\bV'{}^a\op V''{}^a$ for all $a$. If $a\le1-e$ then
$V^a=V=V'\op V''=V'{}^a\op V''{}^a$. (We use that $a\le1-e''$.) If $2-e\le a\le e-1$ we have 
$$V^a=(r'\op0)\i(\bV^a)=(r'\op0)(\bV'{}^a\op V''{}^a)=V'{}^a\op V''{}^a.$$
If $a\ge e$ then $V^a=0=V'{}^a\op V''{}^a$. (We use that $a\ge e''$.) 
 
Finally we assume that $e'\ge e''$ (hence $e'=e\ge2$) and $\l\ne0$ (hence $\l'\ne0$). Then $L=L'\op 0$, 
$\bV=\bV'\op V''$. By the induction hypothesis we have $\bV^a=\bV'{}^a\op V''{}^a$ for all $a$. If $a\le-e$ then 
$V^a=V=V'\op V''=V'{}^a\op V''{}^a$. (We use that $a\le1-e''$.) If $1-e\le a\le e$, we have 
$$V^a=(r'+0)\i(\bV^a)=(r'+0)\i(\bV'{}^a\op V''{}^a)=V'{}^a\op V''{}^a.$$
If $a\ge e+1$ then $V^a=0=V'{}^a\op V''{}^a$. (We use that $a\ge e''$.) 

This completes the proof of (a).

\subhead 3.4\endsubhead
In this subsection we assume that $V,Q,(,),\nn,Is(V)$ are as in 1.2 with $p=2$ and $\nn\ge3$. Let $SO(V)$ be
the identity component of $Is(V)$.
Let $V_0$ be a $\CC$-vector space of dimension $\nn$ with a fixed nondegenerate symmetric bilinear form $(,)$. 
Let $SO(V_0)$ be the corresponding special orthogonal group. Let $U_{V_0}$ be the set of unipotent elements in 
$SO(V_0)$. For any $u_0\in U_{V_0}$ we define subspaces $V_0^a=(V_0)_{u_0}^a$ ($a\in\ZZ$) of $V_0$ in the same 
way as $V^a_u$ were defined in 3.5 for $u\in U_V$ except that we now take $\l$ to be always zero (compare
\cite{\USCII, 3.3}). 

Assume that we are given a direct sum decomposition
$$V=V'\op V''_1\op V''_2\op\do\op V''_{2k-1}\op V''_{2k}$$
where $V',V''_i$ are subspaces of $V$ such that $Q|_{V''_i}=0$, $(V',V''_i)=0$ for $i\in[1,2k]$,
$(V''_i,V''_j)=0$ if $i+j\ne2k+1$. Let 
$$V''=V''_1\op V''_2\op\do\op V''_{2k-1}\op V''_{2k}.$$
Let 
$$V_0=V'_0\op V''_{10}\op V''_{20}\op\do\op V''_{2k-1,0}\op V''_{2k,0}$$
be a direct sum decomposition where $V'_0,V''_{i0}$ are subspaces of $V_0$ such that $(,)$ is $0$ on $V''_{i0}$, 
$(V'_0,V''_{i0})=0$ for $i\in[1,2k]$, $(V''_{i0},V''_{j0})=0$ if $i+j\ne2k+1$. Let 
$$V''_0=V''_{10}\op V''_{20}\op\do\op V''_{2k-1,0}\op V''_{2k,0}.$$
Let $L$ (resp. $L_0$) be the simultaneous stabilizer in $SO(V)$ (resp. $SO(V_0)$) of the subspaces $V',V''_i$ 
(resp. $V'_0,V''_{i0}$), $i\in[1,2k]$.

Let $u\in SO(V)$ be such that $uV'=V'$, $uV''_i=V''_i$ for $i\in[1,2k]$ and such that the $SO(V)$-conjugacy class
$\g$ of $u$ is also an $Is(V)$-conjugacy class. Then $uV''=V''$. Let $u'\in SO(V')$, $u''\in SO(V'')$ be the 
restrictions of $u$. Note that $u''$ is split. Let $\g_1$ be the conjugacy class of $u$ in $L$. Let 
$\g'_1=\ti\r_1(\g_1)$, a unipotent class in $L_0$; here $\ti\r_1:\uuL@>>>\uuL_0$ is defined like $\ti\r$ in 3.1 
but in terms of $L,L_0$ instead of $G,G_0$. Let $u_0\in\g'_1$. Let $u'_0,u''_0$ be the restrictions of $u_0$ to 
$V'_0,V''_0$. By results in \cite{\USCII, 2.9}, we have 
$$\dim V'_{0,u'_0}{}^a=\dim V'_{u'}{}^a,\qua\dim V''_{0,u''_0}{}^a=\dim V''_{u''}{}^a$$
for all $a$.

Let $\g_0$ be the conjugacy class of $u_0$ in $SO(V_0)$. From 3.3(a) and the analogous result for 
$V_0=V'_0\op V''_0$ we see that $\dim V_{0,u_0}{}^a=\dim V_u^a$ for all $a$. From this we deduce, using the 
definitions and the fact that $\g$ is an $Is(V)$-conjugacy class, that $\g_0=\ti\r(\g)$. 

\subhead 3.5\endsubhead
We now prove 3.1(a) assuming that $G=SO(V),G_0=SO(V_0)$ are as in 3.4. Let $C\in\uWW$ be such that $C$ is a 
$W$-conjugacy class. We can find a standard parabolic subgroup $\WW'$ of $\WW$ and an elliptic conjugacy class 
$C'$ of $\WW'$ such that $C'\sub C$. Let $P$ (resp. $P_0$) be a parabolic subgroup of $G$ (resp. $G_0$) of the 
same type as $\WW'$. Let $L$ be a Levi subgroup of $P$ and let $L_0$ be a Levi subgroup of $P_0$. We can assume 
that $L,L_0$ are as in 3.4. Let $\Ph_0^{L_0},\Ph^L$ be the maps analogous to $\Ph_0,\Ph$ with $G_0,G$ replaced by
$L_0,L$. Let $\g_1=\Ph^L(C')$. By definition $\g:=\Ph(C)$ is the unipotent class in $G$ that contains $\g_1$; 
note that $\g$ is an $Is(V)$-conjugacy class. Let $\g'_1=\ti\r_1(\g_1)$ (notation of 3.4), a unipotent class in 
$L_0$. Using 3.2(a) for $L,L_0,C'$ instead of $G,G_0,C$ we see that $\g'_1=\ti\r_1\Ph^L(C')=\Ph_0^{L_0}(C')$. By 
definition, $\g_0:=\Ph_0(C)$ is the unipotent class in $G_0$ that contains $\g'_1$. From 3.4 we have 
$\g_0=\ti\r(\g)$ that is $\Ph_0(C)=\ti\r\Ph(C)$. This completes the proof of 3.1(a) and that of Theorem 0.4.

\subhead 3.6\endsubhead
The equality $\r=\ti\r$ in 3.1 provides an explicit computation of the map $\ti\r$ for special orthogonal groups 
(since the maps $\Ph,\Ps$ are described in each case explicitly in Sections 1,2). The first explicit computation 
of $\ti\r$ in this case was given in \cite{\XUE} in terms of Springer representations instead of the maps 
$\Ph,\Ps$.

\subhead 3.7\endsubhead
According to \cite{\SPA,III,5.4(b)} there is a well defined map $\r':\uuG@>>>\uuG_0,\g\m\g_0$ such that:
$\g\le\p'(\g_0)$ ($\g_0\in\uuG_0$); if $\g\le\p'(\g'_0)$ ($\g'_0\in\uuG_0$) then $\g_0\le\g'_0$. We note the 
following result (see \cite{\XUE}):

(a) {\it If $G$ is almost simple of type $B,C$ or $D$ then $\ti\r=\r'$.}
\nl
Next we note:

(b) {\it For any $G$ we have $\r=\r'$.}
\nl
We can assume that $G$ is almost simple and that $p$ is a bad prime for $G$. If $G$ is of exceptional type, then
$\r$ can be computed from the tables in Section 2 and $\r'$ can be computed from the tables in \cite{\SPA, IV}; 
the result follows. If $G$ is of type $B,C$ or $D$ then as we have seen earlier we have $\ti\r=\r$ and the result
follows from (a).

\subhead 3.8\endsubhead
The results in this subsection are not used elsewhere in this paper.
We define a map $\r'':\uuG@>>>\uuG_0$ as follows. Let $\g\in\uuG$. We can find a Levi subgroup $L$ of a parabolic
subgroup $P$ of $G$ and $\g_1\in\uuL$ such that $\g_1\sub\g$ and $\g_1$ is "distinguished" in $L$ (that is, any 
torus in the centralizer in $L$ of an element in $\g_1$ is contained in the centre of $L$). Let $L_0$ be a Levi 
subgroup of a parabolic subgroup $P_0$ of $G_0$ of the same type as $L$. We have $\g_1=\p_L(\g'_1)$ for a well 
defined $\g'_1\in\uuL_0$ where $\p_L$ is the map analogous to $\p$ in 3.1 but for $L,L_0$ instead of $G,G_0$. Let
$\r''(\g)$ be the unique unipotent class in $G_0$ which contains $\g'_1$. This is independent of the choices and
$\g\m\r''(\g)$ defines the map $\r''$. We show:

(a) $\r''=\r$.
\nl
Let $\g,L,L_0,P,P_0,\g_1,\g'_1$ be as above. Let $\WW'$ be a standard parabolic subgroup of $\WW$ of the same
type as $P,P_0$. We can find an elliptic conjugacy class $C'$ of $\WW'$ such that $\g_1=\Ph^L(C')$,
$\g'_1=\Ph^{L_0}(C')$, where $\Ph^L,\Ph^{L_0}$ are defined like $\Ph$ in terms of $L,L_0$ instead of $G$. Let 
$C$ be the conjugacy class of $\WW$ that contains $\WW'$. From \cite{\WE, 1.1, 4.5} we see that $\Ph(C)$ is the 
unique unipotent class in $G$ that contains $\g_1$ and $\Ph_0(C)$ is the unique unipotent class in $G_0$ that 
contains $\g'_1$. Thus $\Ph(C)=\g$ and $\Ph_0(C)=\r(\g)$. We see that 
$$\r''(\g)=\r''(\Ph(C))=\Ph_0(C)=\r\Ph(C)=\r(\g).$$
This proves (a).

\subhead 3.9\endsubhead
In this subsection we assume that $G$ is of type $B,C$ or $D$ and $p=2$. The identity $\Ps_0=\Ps\p'$ follows from
the explicit combinatorial description of $\Ps_0,\Ps$ given in section 1 and the explicit combinatorial 
description of $\p'$ given in \cite{\SPA, III,6.1,7.2,8.2}. This completes the proof of 3.1(b) hence also that
0.4(b).
 
The fact that $\Ps'_0$ (see 0.2) and $\p'$ might be described by the same combinatorics was noticed by the author
in 1987 who proposed it as a problem to Spaltenstein; he proved it in \cite{\SPAC} (see \cite{\SPAC, p.193}). 
Combining this with the (simple) description of $\Ps$ given in Section 1, we deduce that $\Ps'_0=\Ps\p'$. In 
particular we have $\Ps_0=\Ps'_0$. 

\head 4. Proof of Theorem 0.6\endhead
\subhead 4.1\endsubhead
Let $\uuG^{\sp}$ be the image of $\uuG_0^{\sp}$ (see 0.5) under the imbedding $\p':\uuG_0@>>>\uuG$ (see 3.1). The 
unipotent classes in $\uuG^{\sp}$ are said to be special. The following result can be extracted from 
\cite{\SPA, III}.

(a) {\it There exists order preserving maps $e:\uuG@>>>\uuG$, $e_0:\uuG_0@>>>\uuG_0$ such that the image of $e$ 
(resp. $e_0$) is equal to $\uuG^{\sp}$ (resp. $\uuG_0^{\sp}$) and such that for any $\g\in\uuG$, (resp. 
$\g_0\in\uuG^0$), $\g\le e(\g)$ (resp. $\g_0\le e_0(\g_0)$). The map $e$ (resp. $e_0$) is unique. Moreover we 
have $e^2=e$, $e_0^2=e_0$, $e=\p'e_0\r'$ where $\r'$ is as in 3.7 .}
\nl
Strictly speaking (a) does not appear in \cite{\SPA, III} in the form stated above since the notion of special 
representations from \cite{\LU} and the related notion of special unipotent class do not explicitly appear in 
\cite{\SPA, III} (although they served as a motivation for Spaltenstein, see \cite{\SPA, III,9.4}). Actually (a) 
is a reformulation of results in \cite{\SPA, III} taking into account developments in the theory of Springer 
representations which occured after \cite{\SPA, III} was written.

We now discuss (a) assuming that $G$ is almost simple. Let $d_0:\uuG_0@>>>\uuG_0$ be a map as in 
\cite{\SPA, III,1.4}. If $G$ is of exceptional type we further require that the image of $d_0$ has a minimum 
number of elements (see \cite{\SPA, III,9.4}). Let $e_0=d^2_0$. Then $e_0$ is order preserving and
$\g_0\le e_0(\g_0)$ for any $\g_0\in\uuG_0$. Moreover if we set $e=\p'e_0\r'$ then $e$ is order preserving and 
$\g\le e(\g)$ for any $\g\in\uuG$ (see \cite{\SPA, III,5.6}). If $G$ is of type $B,C$ or $D$ then the map $e_0$ 
is described combinatorially in \cite{\SPA, III,3.10} hence its image is explicitly known; using the explicit 
description of the Springer correspondence given in these cases in \cite{\ICC} we see that this image is exactly 
$\uuG_0^{\sp}$. (See also \cite{\SPA, III,3.11}.) If $G$ is of exceptional type then the image of $e_0$ is 
described explicitly in the tables in \cite{\SPA, p.247-250} and one can again check that it is exactly 
$\uuG_0^{\sp}$. Then the image of $e$ is automatically $\uuG^{\sp}$. The uniqueness in (a) is discussed in 4.2.

\subhead 4.2\endsubhead
Let $\g\in\uuG$. We show that 

(a) {\it there is a unique element $\te(\g)\in\uuG^{\sp}$ such that:

($*$) $\g\le\te(\g)$; if $\g\le\g'$, $(\g'\in\uuG^{\sp})$, then $\te(\g)\le\g'$.
\nl
Moreover we have $\te(\te(\g))=\te(\g)$.}
\nl
We show that $\te(\g)=e(\g)$ satisfies ($*$). We have $\g\le e(\g)$. If $\g'\in\uuG^{\sp}$ and $\g\le\g'$ then
$e(\g)\le e(\g')=\g'$ (since $e^2=e$), as required. Conversely let $\te(\g)$ be as in ($*$). From $\g\le\te(\g)$ 
we deduce $e(\g)\le e(\te(\g))$ hence  $e(\g)\le\te(\g)$. On the other hand taking $\g'=e(g)$ in ($*$) (which 
satisfies $\g\le\g'$) we have $\te(\g)\le e(\g)$ hence $\te(\g)=e(\g)$. This proves (a) and that 

(b) $\te(\g)=e(\g)$.
\nl
Note that (a),(b) and the analogous statements for $G_0$ establish the uniqueness statement in 4.1(a) and the
identity

(c) $\te=\p'\te_0\r'$
\nl
where $\te_0:\uuG_0@>>>\uuG_0$ is the map analogous to $\te:\uuG@>>>\uuG$ (for $G_0$ instead of $G$); we have
$\te_0=e_0$. Note that another proof of (c) which does not rely on the results of \cite{\SPA, III} is given in 
\cite{\XUE}.

According to \cite{\UNI}, for any $\g\in\uuG_0$ we have

(d) $E_{\te_0(\g)}=[E_\g]$;

\subhead 4.3\endsubhead
In this subsection we assume that $p=2$, $G,G_0$ are simple adjoint of type $B_n$ ($n\ge3$) and let $G^*,G^*_0$ 
be almost simple simply connected groups of type $C_n$ defined over $\kk,\CC$ respectively. Note that $G,G^*$ 
have the same Weyl group $\WW$. Define $\x_0:\uuG_0@>>>\cs_\WW$ by $\g\m[E_\g]$ (notation of 0.5). Define 
$\x:\uuG@>>>\cs_\WW$ by $\g\m(\text{ Springer representation attached to }\te(\g))$. We write 

$\uuG^*,\uuG^*_0,\Ph^*,\Ph^*_0,\r^*,\r'{}^*,\x_0^*,\x^*$
\nl
for the analogues of 

$\uuG,\uuG_0,\Ph,\Ph_0,\r,\r',\x_0,\x$
\nl
with $G,G_0$ replaced by $G^*,G^*_0$. We show that

(a) $\x_0\r'(\g)=\x(\g)$
\nl
for any $\g\in\uuG$. The right hand side is the Springer representation attached to $\te(\g)$. The left hand side
is $[E_{\r'(\g)}]$ which by 4.2(d) is equal to $E_{\te_0(\r'(\g))}$; by the definition of $\p'$ this equals the 
Springer representation attached to $\p'\te_0\r'(\g)=\te(\g)$ (see 4.2(c)), proving (a).

Let $\a:G^*@>>>G$ be the standard exceptional isogeny. Let $\a':\WW@>>>\WW$, $\a'':\cs_\WW@>>>\cs_\WW$ be 
the induced bijections. ($\a',\a''$ are the identity in our case.) We consider the diagram
$$\CD
\uWW@>\Ph_0>>\uuG_0@>=>>\uuG_0@>\x_0>>\cs_\WW\\
@A=AA    @A\r AA     @A\r'AA       @A=AA  \\
\uWW@>\Ph>>\uuG@>=>>\uuG@>\x>>\cs_\WW\\
@V\a'VV    @V\a VV     @V\a VV       @V\a''VV\\
\uWW@>\Ph^*>>\uuG^*@>=>>\uuG^*@>\x^*>>\cs_\WW\\
@V=VV    @V\r^*VV     @V\r'{}^*VV       @V=VV\\
\uWW@>\Ph^*_0>>\uuG^*_0@>=>>\uuG^*_0@>\x^*_0>>\cs_\WW
\endCD$$
The top three squares are commutative by 0.4(a), 3.7(b) and (a) (from left to right). The bottom three squares 
are commutative by 0.4(a), 3.7(b) and (a) (from left to right) with $G$ replaced by $G^*$. The middle three 
squares are commutative by the definition of $\a$. We see that the diagram above is commutative. It follows that 
$\x_0\Ph_0=\x\Ph=\x^*\Ph^*=\x_0^*\Ph_0^*$. In particular, $\x_0\Ph_0=\x_0^*\Ph_0^*$. From the definition we have 
$\ti\Ph_0=\x_0\Ph_0$; similarly we have $\ti\Ph^*_0=\x^*_0\Ph^*_0$. Hence we have $\ti\Ph_0=\ti\Ph^*_0$. This 
proves the last sentence in 0.6(a). (Note that in the case where $C\in\uWW$ is elliptic, the equality 
$\ti\Ph_0(C)=\ti\Ph^*_0(C)$ follows also from \cite{\CSM, 3.6}.) 

\subhead 4.4\endsubhead
We now repeat the arguments of 4.3 in the case where $G=G^*,G_0=G^*_0$ are simple of type $F_4$ with $p=2$ or of 
type $G_2$ with $p=3$ or of type $B_2$ with $p=2$ and $\a:G^*@>>>G$ is the standard exceptional isogeny. In these
cases $\a':\WW@>>>\WW$ is the nontrivial involution of the Coxeter group $\WW$ and 
$\a'':\cs_\WW@>>>\cs_\WW$ is induced by $\a'$. We have $\ti\Ph_0=\ti\Ph_0^*$. As in 4.3 we obtain that 
$\ti\Ph_0\a'=\a''\ti\Ph_0$. This, together with 4.3, implies the validity of Theorem 0.6(a).

\subhead 4.5\endsubhead
In the setup of 4.3 we define $\et_0:\cs_\WW@>>>\uuG_0$ by $E_\g\m\g$ where 
$\g\in\uuG_0^{\sp}$ (notation of 0.5). Define 
$\et:\cs_\WW@>>>\uuG$ by $E_\g\m\p'(\g)$ where $\g\in\uuG_0^{\sp}$.
Let $\p'{}^*,\et^*_0,\et^*$ be the analogues of $\p',\et_0,\et$ with $G,G_0$ replaced by $G^*,G^*_0$. We consider
the diagram
$$\CD
\uWW@<\Ps_0<<\uuG_0@<\et_0<<\cs_\WW\\
@V=VV    @V\p'VV            @A=AA  \\
\uWW@<\Ps<<\uuG@<\et<<\cs_\WW\\
@V\a'VV    @V\a VV     @V\a''VV\\
\uWW@<\Ps^*<<\uuG^*@<\et^*<<\cs_\WW\\
@A=AA    @A\p'{}^*AA      @V=VV\\
\uWW@<\Ph^*_0<<\uuG^*_0@<\et^*_0<<\cs_\WW
\endCD$$
The left top square is commutative by 0.4(b) and by the equality $\p=\p'$. Similarly, the left bottom square is 
commutative. The right top square and the right bottom square are commutative by definition. The middle two 
squares are commutative by the definition of $\a$. We see that the diagram above is commutative. It follows that 
$\Ps_0\et_0=\Ps\et=\Ps^*\et^*=\Ps_0^*\et_0^*$. In particular, $\Ps_0\et_0=\Ps_0^*\et_0^*$. 
Hence we have $\ti\Ps_0=\ti\Ps^*_0$. This proves the last sentence in 0.6(b). 

\subhead 4.6\endsubhead
We now repeat the arguments of 4.5 in the case where $G=G^*,G_0=G^*_0$ are simple of type $F_4$ with $p=2$ or of 
type $G_2$ with $p=3$ or of type $B_2$ with $p=2$ and $\a:G^*@>>>G$ is the standard exceptional isogeny. In these
cases $\a':\WW@>>>\WW$ is the nontrivial involution of the Coxeter group $\WW$ and 
$\a'':\cs_\WW@>>>\cs_\WW$ is induced by $\a'$. We have $\ti\Ps_0=\ti\Ps_0^*$. As in 4.5 we obtain that 
$\ti\Ps_0\a''=\a'\ti\Ps_0$. This, together with 4.5, implies the validity of Theorem 0.6(b).

\head 5. Explicit description of the set of special conjugacy classes in $\WW$\endhead
\subhead 5.1\endsubhead
In this section we assume that $\kk=\CC$. We will describe in each case, assuming that $G$ is simple, the set 
$\uWW_{\sp}$ of special conjugacy classes in $\WW$ and the bijection $\t:\uWW_{\sp}@>\si>>\cs_\WW$ induced by 
$\ti\Ph_0$ in 0.5.

We fix $n\ge2$. Let $A$ be the set of all $(x_1\ge x_2\ge\do)$ where $x_i\in\NN$ are zero for large $i$ such that

$x_1+x_2+\do=2n$,

for any $i\ge1$ we have $x_{2i-1}=x_{2i}\mod2$,

for any $i\ge1$ such that  $x_{2i-1},x_{2i}$ are odd we have $x_{2i-1}=x_{2i}$.
\nl
Let $A'$ be the set of all $((y_1\ge y_2\ge\do),(z_1\ge z_2\ge\do))$ where $y_i,z_i\in\NN$ are zero for large $i$
such that 

$(y_1+y_2+\do)+(z_1+z_2+\do)=n$,  

$y_{i+1}\le z_i\le y_i+1$ for $i\ge1$.
\nl
Define $h:A@>>>A'$ by 

$(x_1\ge x_2\ge\do)\m((y_1\ge y_2\ge),(z_1\ge z_2\ge))$
\nl
where 

$y_i=x_{2i}/2, z_i=x_{2i-1}/2$  if $x_{2i-1},x_{2i}$ are even,

$y_i=(x_{2i}-1)/2, z_i=(x_{2i-1}+1)/2$ if $x_{2i-1}=x_{2i}$ are odd.
\nl
Define $h':A'@>>>A$ by 

$((y_1\ge y_2\ge\do),(z_1\ge z_2\ge\do))\m(x_1\ge x_2\ge\do)$
\nl
where 

$x_{2i}=2y_i,x_{2i-1}=2z_i$ if $z_i\le y_i$,

$x_{2i}=2y_i+1,x_{2i-1}=2z_i-1$ if $z_i=y_i+1$.
\nl
Note that $h,h'$  are inverse bijections.

\subhead 5.2\endsubhead
We preserve the setup of 5.1 and we assume that $G$ is simple of type $C_n$. Let $W$ be the group of permutations
of $[1,2n]$ which commute with the involution $\c:i\m 2n-i+1$. We identify $\WW$ with $W$ as in \cite{\WE, 1.4}. 
To $(x_1\ge x_2\ge\do)\in A$ we associate an element of $W=\WW$ which is a product of disjoint cycles with sizes 
given by the nonzero numbers in $x_1,x_2,\do$ where each cycle of even size is $\c$-stable and each cycle of odd 
size is (necessarily) not $\c$-stable. This identifies $A$ with the subset $\uWW_{\sp}$ of 
$\uWW$. We identify in the standard way $\hWW$ with the
set of all $((y_1\ge y_2\ge\do),(z_1\ge z_2\ge\do))$ where $y_i,z_i\in\NN$ are zero for large $i$ such that 
$(y_1+y_2+\do)+(z_1+z_2+\do)=n$. Under this identification $A'$ becomes $\cs_\WW$ (the special representations of
$\WW$). The bijection $\t:\uWW_{\sp}@>\si>>\cs_\WW$ becomes the bijection $h:A@>\si>>A'$ in 5.1.

\subhead 5.3\endsubhead
Now assume that $G$ is simple of type $B_n$ ($n\ge2$). Let $G_1$ be a simple group of type $C_n$ over $\CC$.
Then $\WW$ can be viewed both as the Weyl group of $G$ and that of $G'$. The bijection
sets $\uWW_{\sp},\cs_\WW$ and the bijection $\t:\uWW_{\sp}@>\si>>\cs_\WW$ is the same from the point of view of
$G'$ as from that of $G$ (see 0.6) hence it is described combinatorially as in 5.2.

In the case where $G$ is simple of type $A_n$ ($n\ge1$) then $\uWW_{\sp}=\uWW,\cs_\WW=\hWW$ are all naturally 
parametrized by partitions of $n$ and the bijection $\t$ is the identity map in these parametrizations.

\subhead 5.4\endsubhead
We fix $n\ge2$. Let $C$ be the set of all $((x_1\ge x_2\ge\do),(e_1,e_2,\do))$ where $x_i\in\NN$ are zero for 
large $i$ and $e_i\in\{0,1\}$ are such that

$x_1+x_2+\do=2n$,

for any $i\ge1$ we have $x_{2i-1}=x_{2i}\mod2$, $e_{2i-1}=e_{2i}$,

for any $i\ge1$ such that $x_{2i-1},x_{2i}$ are odd we have $x_{2i-1}=x_{2i},e_{2i-1}=e_{2i}=0$,

for any $i\ge1$ such that  $x_{2i-1},x_{2i}$ are even and $e_{2i-1}=e_{2i}=0$ we have $x_{2i-1}=x_{2i}$,

for any $i\ge1$ such that  $x_{2i}=0$ we have $x_{2i-1}=0$ and $e_{2i-1}=e_{2i}=0$,

for any $i\ge1$ such that $x_{2i}=x_{2i+1}$ are even we have $e_{2i}=e_{2i+1}=0$.
\nl
Let $C'$ be the set of all $((y_1\ge y_2\ge\do),(z_1\ge z_2\ge\do))$ where $y_i,z_i\in\NN$ are zero for large $i$
such that 

$(y_1+y_2+\do)+(z_1+z_2+\do)=n$, 

$y_{i+1}-1\le z_i\le y_i$ for $i\ge1$.
\nl
Define $k:C@>>>C'$ by 

$((x_1\ge x_2\ge\do),(e_1,e_2,\do))\m((y_1\ge y_2\ge\do),(z_1\ge z_2\ge\do))$
\nl
where 

$y_i=(x_{2i}+1)/2, z_i=(x_{2i-1}-1)/2$ if $x_{2i-1}=x_{2i}$ are odd,

$y_i=x_{2i}/2, z_i=x_{2i-1}/2$  if $x_{2i-1}=x_{2i}$ are even and $e_{2i-1}=e_{2i}=0$,

$y_i=(x_{2i}+2)/2, z_i=(x_{2i-1}-2)/2$  if $x_{2i-1},x_{2i}$ are even and $e_{2i-1}=e_{2i}=1$.
\nl
Define $k':C'@>>>C$ by 

$((y_1\ge y_2\ge\do),(z_1\ge z_2\ge\do))@>>>((x_1\ge x_2\ge\do),(e_1,e_2,\do))$
\nl
where

$x_{2i}=2y_i, x_{2i-1}=2z_i, e_{2i}=e_{2i-1}=0$ if $y_i=z_i$,

$x_{2i}=2y_i-1, x_{2i-1}=2z_i+1, e_{2i}=e_{2i-1}=0$ if $y_i=z_i+1$,

$x_{2i}=2y_i-2, x_{2i-1}=2z_i+2, e_{2i}=e_{2i-1}=1$, if $y_i\ge z_i+2$.
\nl
Note that $k,k'$  are inverse bijections.

Let $C_0$ be the set of all $((x_1\ge x_2\ge\do),(e_1,e_2,\do))$ where $x_i\in\NN$ are zero for 
large $i$ and $e_i\in\{0,1\}$ are such that

$x_1+x_2+\do=2n$,

$x_{2i-1}=x_{2i}$ are even, $e_{2i-1}=e_{2i}=0$ for $i\ge1$.
\nl
Note that $C_0\sub C$.

Let $C'_0$ be the set of all $((y_1\ge y_2\ge\do),(z_1\ge z_2\ge\do))$ where $y_i,z_i\in\NN$ are zero for large 
$i$ such that 

$(y_1+y_2+\do)+(z_1+z_2+\do)=n$, 

$y_i=z_i$ for $i\ge1$.
\nl
Note that $C'_0\sub C'$. Now $k:C@>\si>>C'$ restricts to a bijection $k_0:C_0@>\si>>C'_0$. It maps 
$((y_1\ge y_2\ge\do),(y_1\ge y_2\ge\do))$ to $((2y_1,2y_1,2y_2,2y_2,\do),(0,0,0,\do))$.

\subhead 5.5\endsubhead
We preserve the setup of 5.4 and we assume that $n\ge3$ and that $G$ is simple of type $D_n$.
Let $W$ be as in 5.2 and let $W'$ be the subgroup of $W$ consisting of even permutations of $[1,2n]$.
We identify $\WW=W'$ as in \cite{\WE, 1.4, 1.5}. Let $\c:[1,2n]@>>>[1,2n]$ be as in 5.2.
Let $\uWW_{\sp,0}$ (resp. $\uWW_{\sp,1}$) be the set of special conjugay classes in $\WW$ which are not 
conjugacy classes of $W$ (resp. form a single conjugacy class in $W$). 
Let $\cs_{\WW,0}$ (resp. $\cs_{\WW,1}$) be the set of special representations of $\WW$ which do not extend
(resp. extend) to $W$-modules.

To $(x_1\ge x_2\ge\do)\in C$ we associate an element of $W'=\WW$ which is a product of disjoint cycles with sizes 
given by the nonzero numbers in $x_1,x_2,\do$ where each cycle of even size is $\c$-stable and each cycle of odd 
size is (necessarily) not $\c$-stable. This identifies $C-C_0$ with $\uWW_{\sp,1}$ and $C_0$ with $\uWW_{\sp,0}$ 
modulo the fixed point free involution given by conjugation by an element in $W-W'$.

An element $((y_1\ge y_2\ge\do),(z_1\ge z_2\ge\do))\in C'$ can be viewed as in 5.2 as an irreducible
representations of $W$. This identifies $C'-C'_0$ with the $\cs_{\WW,1}$ and $C'_0$ with $\cs_{\WW,0}$ modulo the
fixed point free involution given by $E\m E'$ where $E\op E'$ extends to a $W$-module. Under these 
identifications, the bijection $\t:\uWW_{\sp}@>\si>>\cs_\WW$ becomes the bijection 
$(C-C_0)\sqc C_0\sqc C_0@>\si>>(C'-C'_0)\sqc C'_0\sqc C'_0$, $(a,b,c)\m(k(a),k_0(b),k_0(c))$.

\subhead 5.6\endsubhead
In 5.7-5.11 we describe explicitly the bijection $\t:\uWW_{\sp}@>\si>>\cs_\WW$ (in the case where $G$ simple of 
exceptional type) in the form of a list of data $\a\m\b$ where $\a$ is a special conjugacy class and $\b$ is a
special representation (we use notation  of \cite{\CA} for the conjugacy classes in $\WW$ and the notation of
\cite{\SPA} for the objects of $\hWW$).

\subhead 5.7 Type $G_2$\endsubhead
We have:

$A_0\m \e$

$A_2\m \th'$   

$G_2\m 1$.

\subhead 5.8. Type $F_4$\endsubhead
We have:

$A_0\m\c_{1,4}$

$2A_1\m\c_{4,4}$

$4A_1\m\c_{9,4}$

$A_2\m\c_{8,4}$

$\tA_2\m\c_{8,2}$

$D_4(a_1)\m\c_{12}$

$D_4\m\c_{8,1}$

$C_3+A_1\m\c_{8,3}$

$F_4(a_1)\m\c_{9,1}$

$B_4\m\c_{4,1}$

$F_4\m\c_{1,1}$.

\subhead 5.9. Type $E_6$\endsubhead
We have:

$A_0\m1_{36}$ 

$A_1\m6_{25}$

$2A_1\m20_{20}$ 

$A_2\m30_{15}$

$A_2+A_1\m64_{13}$ 

$2A_2\m24_{12}$

$A_2+2A_1\m 60_{11}$

$A_3\m81_{10}$

$D_4(a_1)\m80_7$ 

$A_4\m81_6$ 

$D_4\m24_6$

$A_4+A_1\m60_5$    

$D_5(a_1)\m64_4$ 

$E_6(a_2)\m30_3$   

$D_5\m20_2$ 

$E_6(a_1)\m6_1$ 

$E_6\m1_0$.

\subhead 5.10. Type $E_7$\endsubhead
We have:

$A_0\m 1_{63}$

$A_1\m 7_{46}$ 

$2A_1\m 27_{37}$

$(3A_1)'\m 21_{36}$ 

$A_2\m 56_{30}$  

$A_2+A_1\m 120_{25}$ 

$A_2+2A_1\m 189_{22}$ 

$A_3\m 210_{21}$ 

$2A_2\m 168_{21}$
 
$A_2+3A_1\m 105_{21}$  

$(A_3+A_1)'\m 189_{20}$ 

$D_4(a_1)\m 315_{16}$ 

$D_4\m 105_{15}$ 

$D_4(a_1)+A_1\m 405_{15}$ 

$D_4(a_1)+2A_1\m 378_{14}$  

$2A_3+A_1\m 210_{13}$

$A_4\m 420_{13}$ 

$A'_5\m 105_{12}$ 

$A_4+A_1\m 512_{11}$ 

$D_5(a_1)\m 420_{10}$  

$A_4+A_2\m 210_{10}$

$D_5(a_1)+A_1\m 378_9$ 

$E_6(a_2)\m 405_8$ 

$E_7(a_4)\m 315_7$ 

$D_5\m 189_7$   

$A_6\m 105_6$

$D_5+A_1\m 168_6$ 

$D_6(a_1)\m 210_6$ 

$A_7\m 189_5$ 

$E_6(a_1)\m 120_4$  

$E_6\m 21_3$ 

$E_7(a_3)\m 56_3$  

$E_7(a_2)\m 27_2$ 

$E_7(a_1)\m 7_1$ 

$E_7\m 1_0$.

\subhead 5.11. Type $E_8$\endsubhead
We have:

$A_0\m 1_{120}$

$A_1\m 8_{91}$

$2A_1\m 35_{74}$

$A_2\m 112_{63}$                             

$A_2+A_1\m 210_{52}$              

$A_2+2A_1\m 560_{47}$

$A_3\m 567_{46}$

$2A_2\m 700_{42}$                     

$D_4(a_1)\m 1400_{37}$

$D_4(a_1)+A_1\m 1400_{32}$

$D_4\m 525_{36}$

$(2A_3)'\m 3240_{31}$

$A_4\m 2268_{30}$

$D_4(a_1)+A_2\m 2240_{28}$

$A_4+A_1\m 4096_{26}$

$A_4+2A_1\m 4200_{24}$

$D_5(a_1)\m 2800_{25}$

$A_4+A_2\m 4536_{23}$

$A_4+A_2+A_1\m 2835_{22}$

$D_5(a_1)+A_1\m 6075_{22}$

$D_4+A_3\m 4200_{21}$

$E_6(a_2)\m 5600_{21}$

$E_8(a_8)\m 4480_{16}$

$D_5\m 2100_{20}$

$D_6(a_1)\m 5600_{15}$

$A_6\m 4200_{15}$

$A'_7\m 6075_{14}$

$A_6+A_1\m  2835_{14}$

$A_7+A_1\m  4536_{13}$

$E_6(a_1)\m 2800_{13}$

$D_7(a_2)\m 4200_{12}$

$E_6(a_1)+A_1\m 4096_{11}$

$E_6\m 525_{12}$

$E_7(a_3)\m 2268_{10}$

$A_8\m 2240_{10}$

$D_8(a_2)\m 3240_9$

$E_8(a_6)\m 1400_8$

$E_8(a_7)\m 1400_7$

$E_8(a_3)\m 700_6$

$E_7(a_1)\m 567_6$

$D_8\m 560_5$

$E_8(a_5)\m  210_4$

$E_8(a_4)\m 112_3$

$E_8(a_2)\m 35_2$

$E_8(a_1)\m 8_1$

$ E_8\m 1_0$.

\subhead 5.12\endsubhead
From the tables above we see that if $\WW$ is of type $G_2,F_4,E_8$ and $<c^k>$ is the conjugacy class of the
$k$-th power of a Coxeter element of $\WW$ then:

(type $G_2$) $<c^1>,<c^2>$ are special, $<c^3>$ is not special;

(type $F_4$) $<c^1>,<c^2>,<c^3>$ are special, $<c^4>$ is not special;

(type $E_8$) $<c^1>,<c^2>,<c^3>,<c^4>,<c^5>$ are special, $<c^6>$ is not special.
\nl
The numbers $2,3,5$ above may be explained by the fact that the Coxeter number of $\WW$ of type $G_2,F_4,E_8$ is 
$2\T3,3\T4,5\T6$ respectively.

\subhead 5.13\endsubhead
Theorem 0.6 suggests that the set $\uWW_{\sp}$ should also make sense in the case where $\WW$ is replaced by a 
finite Coxeter group $\G$ (not necessarily a Weyl group). Assume for example that $\G$ is a dihedral group of 
order $2m$, $m\ge3$, with standard generators $s_1,s_2$. We expect that if $m$ is odd then $\un\G_{\sp}$ consists
of $1$, the conjugacy class of $s_1s_2$ and the conjugacy class containing $s_1$ and $s_2$; if $m$ is even then 
$\un\G_{\sp}$ consists of $1$, the conjugacy class of $s_1s_2$ and the conjugacy class of $s_1s_2s_1s_2$. This 
agrees with the already known cases when $m=3,4,6$.

\enddocument